  \documentclass[12pt, a4paper]{amsart}

\address{ }
\usepackage[left=3cm,top=3.5cm,right=3cm,marginparwidth=2.5cm]{geometry}
\usepackage{parskip, color}
\usepackage{geometry}

\usepackage{cite}

 \usepackage{hyperref}

\usepackage{marginnote}

\usepackage{tikz}

\usepackage{pdfpages}

\usepackage{setspace}

\theoremstyle{remark}
\newtheorem*{definition}{Definition}

\newtheorem{theorem}{Theorem}[section]
\newtheorem{lemma}[theorem]{Lemma}
\newtheorem{prop}[theorem]{Proposition}

\author{Xuan Hien Nguyen}
\address[Xuan Hien Nguyen]{Iowa State University}
\email{\href{mailto: xhnguyen@iastate.edu}{\nolinkurl{xhnguyen@iastate.edu}}}

\author{Alina Stancu}
\address[Alina Stancu]{Concordia University}
\email{\href{mailto: alina.stancu@concordia.ca}{\nolinkurl{alina.stancu@concordia.ca}}}

\author{Guofang Wei}
\address[Guofang Wei]{UC  Santa Barbara}
\email{\href{mailto: wei@math.ucsb.edu}{\nolinkurl{wei@math.ucsb.edu}}}

\begin{document}

\title[Fundamental gap of horoconvex domains in $\mathbb H^n$]{The fundamental gap of horoconvex domains in $\mathbb H^n$ }
\maketitle

\begin{abstract}
  We show that,  for horoconvex domains in the hyperbolic space, the product of their fundamental gap with the square of their diameter has no positive lower bound. The result follows from the study of the fundamental gap of geodesic balls as the radius goes to infinity. In the process, we improve the lower bound for the first eigenvalue of balls in hyperbolic space. 

\end{abstract}

\section{Introduction}

In this article, the fundamental gap of a domain is the difference between the first two eigenvalues of the Laplacian with zero Dirichlet boundary conditions.   For convex domains in $\mathbb R^n$ or $\mathbb S^n$, $n\geq 2$, it is known from \cite{fundamental, seto2019sharp,dai2018fundamental,he2017fundamental} that $\lambda_2 - \lambda_1 \ge 3\pi^2/D^2$, where $D$ is the diameter of the domain. 

In hyperbolic space, this quantity behaves very differently from the Euclidean and spherical cases.  Recently, the authors showed \cite{BCNSWW2} that for any fixed $D>0$, there are convex domains with diameter $D$ in $\mathbb H^n$, $n\geq 2$, such that  $D^2(\lambda_2 - \lambda_1)$ is arbitrarily small. Since convexity does not provide a lower bound, one naturally asks if imposing a stronger notion of convexity, such as horoconvexity, would imply an estimate for $D^2(\lambda_2 - \lambda_1)$ from below. Recall that for a domain with smooth boundary, convexity corresponds to nonnegative principal curvatures of the boundary, while horoconvexity corresponds to principal curvatures greater or equal to 1.  We show that the quantity $D^2 (\lambda_2 - \lambda_1)$ still tends to zero for all horoconvex domains in hyperbolic space when the diameter tends to infinity. 

 \begin{theorem}  \label{gap-horoconvex}
	For every $n \geq 2$, there exists a constant $C(n)$ such that the Dirichlet fundamental gap of every horoconvex domain $\Omega$ with diameter $D \geq 4 \ln 2$ satisfies
	\[
	\lambda_2(\Omega)- \lambda_1(\Omega) \leq \frac{C(n)}{D^3}. 
	\]
		In particular, as $D \to \infty$, the quantity $(\lambda_2 - \lambda_1) D^2$ tends to $0$. 
\end{theorem}

We prove this by first obtaining the following estimate for the fundamental gap for special horoconvex domains, the 
geodesic balls  in hyperbolic space.

\begin{theorem} \label{thm-main}
	Let $B_R$ be the geodesic ball of radius $R$ in $\mathbb H^n$ and $\lambda_i(B_R)$ be the $i$-th eigenvalue of the Laplace operator $-\Delta$ in $B_R$ with Dirichlet boundary conditions. Then there is a constant $C(n)$ so that 
	\begin{equation}
	\lambda_2 (B_R) - \lambda_1(B_R) \leq \frac{C(n)}{R^3}. \label{gap-l-u}
	\end{equation}
	In particular, as $R \to \infty$, the quantity $(\lambda_2 - \lambda_1) R^2$ tends to $0$. 
\end{theorem}

In the authors' earlier work \cite{BCNSWW2}, it was shown that,  for any fixed $D>0$, one can find a domain $\Omega$ for which $(\lambda_2 (\Omega) - \lambda_1 (\Omega)) D^2$ can be made arbitrarily small. The domains $\Omega \subset \mathbb H^n$ in \cite{BCNSWW2} are convex, but not horoconvex. Their first eigenfunction is not log-concave either.  In contrast, note that the first eigenfunction of $B_R$ is log-concave (see \cite[Corollary 1.1]{IST} and Lemma~\ref{lem: log}). On the one hand, while the log-concavity of the first eigenfunction plays a very important role in estimating the fundamental gap of convex domains in the Euclidean space and sphere, Theorem~\ref{thm-main} shows that the log-concavity of the first eigenfunction in the hyperbolic case does not imply a lower bound estimate for $(\lambda_2 - \lambda_1) D^2$. 
On the other hand, we believe that $D^2$ is not the appropriate factor for domains in the hyperbolic space and we conjecture that, 
for all horoconvex convex domains $\Omega \subset \mathbb H^n$,  we have $\lambda_2 (\Omega) - \lambda_1 (\Omega) \ge c(n, D)$ for some function $c(n, D)$ depending on the  dimension and diameter, that can lead to a lower bound on the fundamental gap appropriately compared with the diameter. This is  true for balls in $\mathbb H^n$, see \eqref{ball-gap-lower}.

Theorem~\ref{thm-main} is proved by transforming the eigenvalue equation of balls to the eigenvalue equation of a Schr\"odinger operator. As a result, we  obtain some immediate upper and lower bound estimates on the first two eigenvalues of balls, which improve and simplify earlier estimates on the first eigenvalues of balls. See Sections 2, 3.

To prove Theorem~\ref{gap-horoconvex}, we exploit the fact  that all big horoconvex domains contain a large ball \cite{Borisenko-Miquel}, see Theorem~\ref{horoconvex-ball}. We then combine Theorem~\ref{thm-main} with
 Benguria and Linde's \cite{Benguria-Linde2007} comparison result for the
 fundamental gap to conclude the proof, see Section 4.

\section{Basic Facts on Eigenvalues of Balls in $\mathbb H^n$}

Here we review some basic facts about first two Dirichlet eigenvalues of balls in the hyperbolic space.  By transforming the eigenvalue equation of balls to its  Schr\"odinger form,  we  obtain some immediate upper and lower bound estimates on the first two eigenvalues  which improve and simplify earlier estimates. 

\subsection{The first eigenvalue}


In this section, let $\lambda_i$ be  the $i$-th eigenvalue of the Laplacian, with Dirichlet boundary conditions, of geodesic balls with radius $r$ in $\mathbb H^n$. 

By \cite{Chavel, Benguria-Linde2007}, the first eigenvalue $\lambda_1$ is the 
first eigenvalue of the $1$-dimensional problem on $[0, r]$
\begin{equation} u'' + \frac{n-1}{\tanh t} u' + \lambda u =0, \  \ u(r) =0, \ u'(0) =0. \label{1-eigenvalue-ball}
\end{equation} 
With the change of variable $u(t) = (\sinh t)^{\frac{1-n}{2}} \bar u(t)$, we have the associated Schr\"odinger equation
\begin{equation}
\label{1-Sch-ball}
-\frac{d^2}{dt^2}\bar u + \frac{n-1}{4} \left( n-1 +\frac{n-3}{\sinh^2 t} \right)  \bar u= \lambda \bar u
\end{equation}
with Dirichlet boundary conditions at $0$ and $r$, and $\lambda_1$ is the  first eigenvalue of \eqref{1-Sch-ball}. Note that the nonconstant potential term changes sign at $n =3$. We immediately notice that, when $n=3$, $\lambda_1 = 1+ \frac{\pi^2}{r^2}$. Since $\sinh^{-2} t \ge \sinh^{-2} r$ on $(0, r]$, the ODE comparison theorem implies:
\begin{lemma} For  $n>3$, 
	\begin{equation}
	\label{lambda1-lower-bound-4}
	\lambda_1 >  \frac{(n-1)^2}{4}+ \frac{\pi^2}{r^2} +  \frac{(n-1)(n-3)}{4\sinh^2 r}.
	\end{equation}
	For $n=2$, 
	\[
	\lambda_1 \le \frac 14+ \frac{\pi^2}{r^2} - \frac1{4 \sinh^{2} r}.  \]
\end{lemma}

The lower bound is sharper than the estimate of \cite[(1.7)]{Artamoshin2016}, which followed the earlier estimate of McKean \cite{McKean1970}. It is also an improvement over \cite[Theorem 5.6]{Savo2009} and an earlier estimate in  \cite[Theorem 5.2]{Gage1980} when $r$ is large and $n>3$. The upper bound in the case $n=2$ is that found by Gage  \cite[Theorem 5.2]{Gage1980}.


%
The bounds in the other direction do not follow directly from the Schr\"odinger equation \eqref{1-Sch-ball}.  
In \cite[Theorem 5.6]{Savo2009}  the following uniform upper and lower bounds for the first eigenvalue $\lambda_1$ is obtained for all $n \ge 2$: 
\begin{equation}
\frac{(n-1)^2}{4}+ \frac{\pi^2}{r^2} - \frac{4\pi^2}{(n-1)r^3} \le \lambda_1 \le  \frac{(n-1)^2}{4}+ \frac{\pi^2}{r^2} + \frac{C}{r^3},  \label{lambda1-all-n}
\end{equation}
with $C = \frac{\pi^2 (n^2-1)}{2} \int_0^\infty \frac{t^2}{\sinh^2 t} dt =  \frac{\pi^4 (n^2-1)}{12}$. 

We will use this lower bound and improve the upper bound in Section 3. 
\subsection{The second eigenvalue}

 The second eigenvalue $\lambda_2$ is studied in \cite[Lemma 3.1]{Benguria-Linde2007}, where it is shown that it is the first eigenvalue of the following equation (see also \eqref{eigenvalue-ball} with $k=1, l=1$): 
\begin{equation}  u'' + \frac{n-1}{\tanh t} u' - \frac{n-1}{\sinh^2 t} u + \lambda u =0, \ \  \ u(r) =0, \ u (t) \sim t \ \mbox{as}\ t \rightarrow 0.  \label{2-eigenvalue-ball}
\end{equation}

Again with the change of variable $u(t) = (\sinh t)^{\frac{1-n}{2}} \bar u(t)$, we have the associated Schr\"odinger equation
\begin{equation}
\label{2-Sch-ball}
-\frac{d^2}{dt^2}  \bar u+ \frac{n-1}{4} \left( n-1 +\frac{n+1}{\sinh^2 t} \right)   \bar u= \lambda \bar u
\end{equation}
 with Dirichlet boundary conditions at $0$ and $r$, where  the second eigenvalue $\lambda_2$ is the first eigenvalue of \eqref{2-Sch-ball}. Using once more the ODE comparison theorem, we obtain
 \begin{equation}
 \lambda_2 \ge
    \frac{(n-1)^2}{4}+ \frac{\pi^2}{r^2} + \frac{n^2-1}{4\sinh^2 r}. \label{lambda2-lower}
    \end{equation}
 To find an upper bound estimate for $\lambda_2$, we will seek in the next section an upper bound for the first eigenvalue of a more general Schr\"odinger equation and, as such, we will simultaneously obtain an upper bound for $\lambda_1$,  slightly improve the one in \eqref{lambda1-all-n}.
 
 From \eqref{1-Sch-ball} and \eqref{2-Sch-ball} we immediately have the following lower bound on the fundamental gap of the ball $B_R \subset \mathbb H^n$ for all $n \ge 2$.
 \begin{equation}  \label{ball-gap-lower}
 	\lambda_2 - \lambda_1 \ge \frac{n-1}{\sinh^2 R}.\end{equation}


 


%
      \section{First Eigenvalue Upper Bound for Schr\"odinger Equation}  \label{sec3}
 
 Let $\lambda_1^\alpha$ be the first eigenvalue of the following equation    
       \begin{equation}
      -\frac{d^2}{dt^2}   u+ \frac{n-1}{4} \left( n-1 +\frac{\alpha}{\sinh^2 t} \right)   u= \lambda  u  \label{Sch-alpha}
      \end{equation}
      with Dirichlet boundary conditions at $0$ and $r$.
  \begin{prop} \label{prop-1}
  	For $\alpha \ge 0$,  we have
  	\begin{equation}
  	\lambda_1^\alpha <  \frac{(n-1)^2}{4}+ \frac{\pi^2}{r^2} + \frac{(n-1) \alpha}{12 r^3} \pi^4.  \label{lambda1-alpha}
  	\end{equation}
  	In particular, the first two eigenvalues of the geodesic ball of radius $r$ in $\mathbb{H}^n$ satisfy
  	\begin{eqnarray}
  		\lambda_1 & < &  \frac{(n-1)^2}{4}+ \frac{\pi^2}{r^2} +  \frac{(n-1) (n-3)}{12r^3} \pi^4, \ \mbox{for} \ n \ge 3, \label{lambda1-ub}\\ 
  	\lambda_2 & <  & \frac{(n-1)^2}{4}+ \frac{\pi^2}{r^2} +  \frac{(n-1) (n+1)}{12r^3} \pi^4, \ \mbox{for} \ n \ge 2.  \label{lambda2-upper}
  	\end{eqnarray}
  \end{prop}         
The upper bound \eqref{lambda1-ub} improves the upper bound in \cite[Theorem 5.6]{Savo2009}, see \eqref{lambda1-all-n}. 

\begin{proof}
	The first Dirichlet eigenvalue of a Schr\"odinger operator $-u''+ Vu$ is a minimizer of the Rayleigh quotient 
	     \[ R[u]=\frac{\int |u'|^2+ V u^2  }{\int u^2}, \]
	     among all non-constant $u$ with $u(0) = u(r)=0$. 
	     
      The equation \eqref{Sch-alpha} with $\alpha =0$ has  its first eigenfunction equal to $v = \sqrt{\frac{2}{r}} \sin (\pi t/r)$.  It is normalized so that $\int_0^r v^2 dt = 1$. 
      Therefore by inserting $v$ into the  Rayleigh quotient associated to \eqref{Sch-alpha}, 
      we find
      \begin{align*}
      \lambda_1^\alpha &\leq \frac{(n-1)^2}{4}+ \int_0^r \left(  \frac{dv}{dt}\right)^2 dt +\int_0^r   \frac{(n-1)\alpha}{4(\sinh t)^2}   v^2 \,dt \\
      &= \frac{(n-1)^2}{4}+ \frac{\pi^2}{r^2} +   \frac{(n-1) \alpha}{4}  \int_0^r   \frac{v^2}{(\sinh t)^2}    \,dt.
      \end{align*}
      
      Using $\sin |x| \leq |x|$, we have 
      	\[
      	r^2 \int_0^r \left( \frac{\sin\left( \pi t/r \right)}{\sinh t} \right)^2 dt \leq \pi^2 \int_0^r \left( \frac{t}{\sinh t} \right)^2 dt < \pi^2 \int_0^\infty \left( \frac{t}{\sinh t} \right)^2 dt = \frac{ \pi^4}{6}.
      	\]
      	This gives $\int_0^r   \frac{v^2}{(\sinh t)^2}    \,dt < \frac{\pi^4}{3r^3}$, hence \eqref{lambda1-alpha}.
      \end{proof}

Combining the lower bound in \eqref{lambda1-all-n} with \eqref{lambda2-upper} gives the estimate \eqref{gap-l-u} in Theorem~\ref{thm-main}.

\section{Horoconvex domains in $\mathbb H^n$}

 A stronger definition of convexity in the hyperbolic space considers horospheres as  natural analogues of Euclidean hyperplanes supporting a convex domain: 
 
  \begin{definition}   	 A set $\Omega \subset \mathbb H^n$ is called horoconvex if,  for every point $p \in \partial \Omega$, there exists a horosphere ${\mathcal{H}}$ 
  through $p$ such that $\Omega$ lies in the horoball bounded by ${\mathcal{H}}$.
  \end{definition}

Recall that a \emph{horosphere} is a sphere with center on the ideal boundary of $\mathbb H^n$ and that a \emph{horoball} is a domain whose boundary is a horosphere. 

When   $\Omega$  is a compact domain with smooth boundary in the hyperbolic space of constant negative curvature $-1$, the domain $\Omega$  is horoconvex if and only if all principal curvatures of the boundary hypersurface  are greater or equal to one. 
As a special case, $B_R$, the geodesic sphere of radius $R$, is horoconvex as each of the principal curvatures of its boundary is equal to $\coth R$, and  $\coth R > 1$ for all $R>0$.  

Finally,  for any compact domain, recall that its inradius is the radius of the largest ball contained in the domain, and that its circumradius is the radius of the smallest ball containing the domain. Part of a result of Borisenko-Miquel \cite[Theorem 1]{Borisenko-Miquel} states the following:

\begin{theorem} \cite{Borisenko-Miquel}
  \label{horoconvex-ball}
Let $\Omega$ be a compact horoconvex domain in $\mathbb{H}^n$ with inradius $r$ and circumradius $R$. Denoting $\tau = \tanh \frac{r}{2}$, then
\begin{equation} \label{eq:BM}
R-r \leq \ln \frac{(1+ \sqrt{\tau})^2}{1 + \tau} < \ln 2, 
\end{equation}  and this bound is sharp.
  \end{theorem}

An immediate consequence of (\ref{eq:BM}) is that the diameter of the domain satisfies $D\le 2R \le 2r+2\ln 2$.  We are now ready to prove Theorem~\ref{gap-horoconvex}.

  \begin{proof} Let $\Omega \subset \mathbb H^n$ be a horoconvex domain of diameter $D$. Choose $R_\Omega$ such that the ball of radius $R_\Omega$ satisfies $\lambda_1(B_{R_\Omega}) = \lambda_1(\Omega)$. Theorem~\ref{horoconvex-ball} implies that $\Omega$ contains a ball of radius $r$ with $r \geq \frac{D}{2} - \ln 2$. By domain monotonicity of the first eigenvalue, 
    $ R_\Omega \ge \frac{D}{2} - \ln 2$, 
 hence   \begin{equation} R_\Omega \ge \frac{D}{4},  \label{R-Omega}
  \end{equation}
  when  $D \ge 4 \ln 2$. 
  	
Using \cite{Benguria-Linde2007}, Benguria-Linde's  upper bound on the second eigenvalue, we have that $$\lambda_2(\Omega) - \lambda_1(\Omega) \leq \lambda_2(B_{R_\Omega}) - \lambda_1(B_{R_\Omega}).$$
    Applying the estimates \eqref{gap-l-u} and \eqref{R-Omega} concludes the proof of Theorem~\ref{gap-horoconvex}.  
  \end{proof}

\section*{Appendix}{{\bf{Small balls and log-concavity of eigenfunction of geodesics balls in $\mathbb M^n_K$}}}
\

To round up the discussion on the fundamental gap of balls in the hyperbolic space, we thought to  include here an observation on the fundamental gap of balls of small radii, as well as a simple argument proving the log-concavity of the first eigenfunction of geodesic balls in simply connected Riemannian manifolds with constant negative sectional curvature.

\subsection{The gap of small balls in negatively curved manifolds}
Let $\mathbb M^n_K$ be the simply connected Riemannian manifold with constant sectional curvature $K$. Here, we assume that  $K$ is negative and write $K =- k^2,\, (k>0)$. Denote by $\lambda_i(n,k,r)$ the eigenvalues of the Laplacian for geodesic balls with radius $r$ in $\mathbb M^n_K$ with Dirichlet boundary condition.

By separation of variables, see \cite{Chavel, Benguria-Linde2007}, the eigenvalues  $\lambda_i(n,k,r)$ are eigenvalues of
\begin{equation}  u'' + \frac{(n-1)k}{ \tanh (kt)} u' - \frac{l(l+n-2)k^2}{\sinh^2 (k t)} u + \lambda u =0,   \label{eigenvalue-ball}
\end{equation}
where $l = 0, 1,2, \cdots$, with boundary condition $u'(0) =0$ for $l =0$, $u(t) \sim t^l$ as $t \to 0$ for $l>0$, and $u(r) =0$.

By scaling, this immediately gives \cite[Lemma 4.1]{Benguria-Linde2007}, for $c>0$,
\begin{equation}
\lambda_i (n, \frac 1c k, cr ) = c^{-2}\lambda_i (n, k, r).
\end{equation}
Hence
\begin{equation}
\lambda_i (n, 1, r ) = r^{-2}\lambda_i (n, r, 1).
\end{equation}
Therefore, for small balls in $\mathbb H^n$, the value $r^2 \lambda_i (n, 1, r )$ is close to the corresponding one in the Euclidean space, as one would expect. Namely,
\begin{lemma}  \label{gap-small-ball}
$$\lim_{r \to 0} r^2 \lambda_i (n, 1, r ) = \lambda_i (n, 0,1) = r^2 \lambda_i (n, 0,r),$$ and 
\begin{equation}
\lim_{r \to 0} r^2 \left(\lambda_2(n,1, r) - \lambda_2(n,1, r) \right)  = r^2(\lambda_2 (n, 0,r)- \lambda_1 (n, 0,r)) = j_{\frac n2, 1}^2 - j_{\frac n2 -1, 1}^2,
\end{equation}
where $j_{p,k}$ is the $k$-th positive zero of the Bessel function $J_p(x)$.
\end{lemma}

\subsection{The first eigenfunction for balls}
The first eigenfunction of balls is purely radial, so it is straightforward to show that it is log-concave, as in 
 the Euclidean and spherical case. 
\begin{lemma} \label{lem: log}
	The first eigenfunction $u_1$ of \eqref{1-eigenvalue-ball} is strictly log-concave.
\end{lemma}
This is in \cite[Corollary 1.1]{IST}, where more general elliptic equations with  power are considered. For convenience, we give a  simple and direct proof here.
\begin{proof}
	First we show $u_1$ is strictly decreasing. Multiplying both sides of  \eqref{1-eigenvalue-ball}  by $\sinh^{n-1} t$, we have
	\[ (u_1' \sinh^{n-1}t)' = -\lambda_1 u_1 \sinh^{n-1}t <0.
	\]
	Since $u_1'(0) =0$, we have $u_1'(t) <0$ for $t \in (0, r)$.
	
	Let $\varphi = (\log u_1)'$. Then $\varphi(0) =0$,  $\varphi <0$ on $(0, r)$, and \[
	\varphi' = \frac{u_1''}{u_1} - \left( \frac{u_1'}{u_1}\right)^2 = -\frac{n-1}{\tanh t} \varphi -\lambda_1 -\varphi^2.   \]
Taking the limit as $t \to 0$ gives $\varphi'(0) = -\lambda_1 - (n-1) \lim_{t \to 0} \frac{\varphi}{\tanh t} = -\lambda_1 - (n-1) \varphi'(0)$. 	Hence, $\varphi'(0) <0$. Now, we claim that $\varphi' (t) < 0$ on $[0, r)$. Otherwise, there exists  $t_1 \in (0, r)$ such that $\varphi' <0$ on $[0, t_1)$, $\varphi' (t_1) = 0$ and $\varphi'' (t_1) \ge 0$. Note that
$\varphi''$ satisfies
\[ \varphi'' = \frac{n-1}{\sinh^2 t} \varphi  -\frac{n-1}{\tanh t} \varphi' - 2\varphi \varphi'.
\]
Evaluating the two sides of the equation at $t_1$ gives \[
0 \le \varphi'' (t_1) = \frac{n-1}{\sinh^2 t_1} \varphi (t_1) <0.  \]
This is a contradiction.	
	\end{proof}

\bibliographystyle{plain}
\bibliography{../VanishingFundamentalGap/references}
\end{document}